\title{Primitive permutation groups containing a cycle}

\author{Gareth A. Jones\\
School of Mathematics\\
University of Southampton\\
Southampton SO17  1BJ, U.K.\\
{\tt G.A.Jones@maths.soton.ac.uk}
}

\documentclass[12pt]{article}
\usepackage[latin1]{inputenc}
\usepackage{a4wide}
\usepackage{latexsym}
\usepackage{amsmath}
\usepackage{amssymb}
\usepackage{amsfonts}
\usepackage{color}
\newtheorem{thm}{Theorem}[section]

\newtheorem{cor}[thm]{Corollary}

\begin{document} 

\maketitle

%\date{}

\begin{abstract}
\noindent The primitive finite permutation groups containing a cycle are classified. Of these, only the alternating and symmetric groups contain a cycle fixing at least three points. The contributions of Jordan and Marggraff to this topic are briefly discussed.
\end{abstract}

%MSC: primary 20B15; secondary 01A55, 20-03, 20B20.

% MSC: primary 20B15 primitive permutation groups; secondary 20B20 multiply transitive groups, 01A55 history of 19th century mathematics, 20--03 history of group theory.

% Key words" primitive group, cycle, fixed points, alternating group.

\section{Introduction}

There is a long tradition, going back to Jordan, of proving that a primitive permutation group of degree $n$, containing an element with a specific cycle structure, must contain $A_n$. The following theorem of Jordan (see~\cite[Theorem~3.3E]{DM} or~\cite[Theorem~13.9]{Wie}) is typical:

\begin{thm} 
Let $G$ be a primitive permutation group of finite degree $n$, containing a cycle of prime length fixing at least three points. Then $G\ge A_n$.
\end{thm}

\noindent The following extension of this result, removing the primality condition, was motivated by a question raised by Alexander Zvonkin in connection with his work with Fedor Pakovich on polynomials and weighted plane trees:

\begin{thm}
Let $G$ be a primitive permutation group of finite degree $n$, not containing the alternating group $A_n$. Suppose that $G$ contains a cycle fixing $k$ points, where $0\le k\le n-2$. Then one of the following holds:
\begin{enumerate}
\item $k=0$ and either
\begin{description}
\item {\rm (a)} $C_p\le G\le AGL_1(p)$ with $n=p$ prime, or
\item {\rm (b)} $PGL_d(q)\le G\le P\Gamma L_d(q)$ with $n=(q^d-1)/(q-1)$ and $d\ge 2$ for some prime power $q$, or
\item {\rm (c)} $G=L_2(11)$, $M_{11}$ or $M_{23}$ with $n=11, 11$ or $23$ respectively.
\end{description}
\item $k=1$ and either
\begin{description}
\item {\rm (a)} $AGL_d(q)\le G\le A\Gamma L_d(q)$ with $n=q^d$ and $d\ge 1$ for some prime power $q$, or
\item {\rm (b)} $G=L_2(p)$ or $PGL_2(p)$ with $n=p+1$ for some prime $p\ge 5$, or
\item {\rm (c)} $G=M_{11}$, $M_{12}$ or $M_{24}$ with $n=12, 12$ or $24$ respectively.
\end{description}
\item $k=2$ and $PGL_2(q)\le G\le P\Gamma L_2(q)$ with $n=q+1$ for some prime power $q$.
\end{enumerate}
\end{thm}

\begin{cor}
Let $G$ be a primitive permutation group of finite degree $n$, containing a cycle with $k$ fixed points. Then $G\ge A_n$ if $k\ge 3$, or if $k=0, 1$ or $2$ and $n$ avoids the values listed in parts (1), (2) or (3) of Theorem~1.2.
\end{cor}

\medskip

\noindent{\bf Comments 1} If a permutation $g$ has a cycle of length coprime to all its other cycle lengths, then some power of $g$ is a cycle of the same length, so these results can be applied to it.

\medskip

\noindent{\bf 2} It is straightforward to check that the groups $G$ listed in Theorem~1.2 all have elements with the appropriate cycle structures. Moreover, they are all primitive. In fact, apart from proper subgroups of $AGL_1(p)$ in 1(a), they are all doubly transitive.

\medskip

\noindent{\bf 3} In general, one cannot remove the hypothesis that $G$ is primitive. For instance, if $m$ is a proper divisor of $n$ then the imprimitive group $S_m\wr S_{n/m}$ of degree $n$ contains a cycle $g$ with $k$ fixed points for $k=m, 2m, \ldots, n-2m$ (permuting the blocks non-trivially), and for $n-m\le k\le n-2$ (leaving each block invariant). However, if $k$ is coprime to $n$ and less than $n/2$, then any transitive group containing $g$ is primitive, so these results apply.

\medskip

\noindent{\bf 4}  A similar result to Theorem~1.2, restricted to the case where $n$ is prime and $k>0$, has been obtained by Bouw and Osserman~\cite[Prop.~3.1]{BO}, who apply it to covers of curves in positive characteristic. (Their proof can be simplified by using the fact that the transitive groups of prime degree $n$ are known: apart from $S_n$ and $A_n$, they are the groups $G$ in parts 1(a) or 1(c) of Theorem~1.2, together with those in 1(b) for which $n$ is prime.)

\medskip

\noindent{\bf 5} As in the preceding comment, some of the motivation for results of this type comes from covering space theory. The monodromy group of a covering is the group of permutations of the sheets obtained by lifting closed paths. It is primitive if and only if the covering is not a composition of coverings of smaller degrees. Local branching information provides cycle structures for certain elements of this group, so it is useful to know which primitive groups contain elements with given cycle structures.

\medskip

\noindent{\bf 6} There is a similar situation in Galois theory: if $f(t)\in{\mathbb Z}[t]$ then the Galois group $G$ of $f$ acts primitively on the roots if and only if $f$ is not a composition of polynomials of smaller degrees. If $p$ is a prime not dividing the discriminant of $f$, the degrees of the irreducible factors of the reduction mod~$(p)$ of $f$ give the cycle structure of an element of $G$.

\medskip

\noindent{\bf 7} Finally, it should be emphasised that the proof of Theorem~1.2 relies heavily on the classification of finite simple groups --- in particular, on the resulting classification of doubly transitive groups. It seems hopeless to expect proofs of results such as this using only the methods available to Jordan and his contemporaries.

\section{Proof of Theorem 1.2}

The case $k=0$ has been dealt with by the author in~\cite{Jon}, completing work of Feit~\cite{Fei}, while the case $k=1$ has been dealt with by M\"uller in~\cite[Theorem~6.2]{Mue96} (see also~\cite[Theorem~3.2]{Mue01} and~\cite{BP}). Thus we may assume that $k\ge 2$.

A theorem of Jordan, often attributed to Marggraff (see Section~3) shows that $G$, being primitive and containing a cycle with $k$ fixed points, is $(k+1)$-transitive; since $k\ge 2$, $G$ is at least $3$-transitive. As a result of the classification of finite simple groups, the multiply transitive finite permutation groups are known (see~\cite{Cam}, for instance). In particular, the $3$-transitive groups $G\not\ge A_n$ are as follows:
\begin{description}
\item (i) various groups $G$ such that $L_2(q)\le G\le P\Gamma L_2(q)$, with $n=q+1$ for some prime power $q$;
\item (ii) various subgroups $G\le AGL_d(2)$ with $n=2^d$ and $d\ge 3$;
\item (iii) $M_{11}$ with $n=11$ or $12$, $M_{12}$ with $n=12$, $M_{22}$ and ${\rm Aut}\,M_{22}$ with $n=22$, $M_{23}$ with $n=23$, $M_{24}$ with $n=24$.
\end{description}
All these groups appear in their natural representations, apart from $M_{11}$ acting on the $n=12$ cosets of a subgroup $L_2(11)$. Of these groups, only $M_{11}, M_{12}, M_{23}$ and $M_{24}$ in their natural representations are $4$-transitive, only $M_{12}$ and $M_{24}$ are $5$-transitive, and none are $6$-transitive. Thus $2\le k\le 4$.

The groups in (iii) can be eliminated since inspection of the groups or of their character tables in~\cite{ATLAS} shows that they do not contain $(n-k)$-cycles for such values of $k$.

If $G\le AGL_d(2)$, as in (ii), then $G$ is only $3$-transitive, so $k=2$. Then the subgroup $G_0$ fixing $0$ is a $2$-transitive subgroup  of $GL_d(2)$ of degree $m=2^d-1$ on the non-zero points, containing an $(m-1)$-cycle, so it must appear in case (2) of this theorem (with $m$ replacing $n$). However, the groups in parts (b) and (c) of (2) all have even degrees, so $G_0$ is an affine group of dimension $c$ and degree $q^c$ for some prime power $q=p^e$, as in (a), with $2^d-1=q^c=p^{ec}$. At this point it is tempting to quote Mih\u ailescu's proof of Catalan's conjecture~\cite{Mih}, but a much simpler argument suffices. Since $d\ge 2$ we have $p^{ec}\equiv -1$ mod~$(4)$, so $ec$ is odd. Thus $p^{ec}+1$ has an odd factor $p^{ec-1}-p^{ec-2}+\cdots+1$, which must be $1$, so $e=c=1$, $q$ is prime and $G_0\cong AGL_1(q)$. This group has a normal Sylow $q$-subgroup $Q\cong C_q$, which permutes the non-zero points transitively, so it is a Singer subgroup of $GL_d(2)$, and hence of index $d$ in its normaliser~\cite[II.7.3]{Hup}. However, $Q$ is normal and of index $q-1$ in $G_0$, so $d$ is divisible by $q-1=2^d-2$, which is impossible for $d\ge 3$.

There remain the groups $G$ in (i). Again, these are only $3$-transitive, so $k=2$. If $q=2^e$ then $L_2(q)=PGL_2(q)$ is $3$-transitive, with an element fixing $0$ and $\infty$ and inducing a cyclic permutation of the remaining $q-1$ points, so the same holds for each $G$ (one for each divisor of $e$). Now suppose that $q=p^e$ for some odd prime $p$. In this case, $P\Gamma L_2(q)/L_2(q)\cong C_2\times C_e$, with the images of $PGL_2(q)$ and $P\Sigma L_2(q)$ giving the direct factors. Now $P\Sigma L_2(q)$ is not $3$-transitive: for instance, the subgroup fixing $0$ and $\infty$ has two orbits on the remaining points, consisting of the quadratic residues and the non-residues in the field ${\mathbb F}_q$. Each element of $P\Gamma L_2(q)\setminus P\Sigma L_2(q)$ transposes the corresponding two orbits of $P\Sigma L_2(q)$ on distinct ordered triples (which are also the orbits of $L_2(q)$), so $G$ is $3$-transitive if and only if it is not a subgroup of $P\Sigma L_2(q)$. It follows that for odd $e$, the $3$-transitive groups $G$ are again those containing $PGL_2(q)$, and each of these contains a cycle with two fixed points. If $e$ is even, however, there is another complement for the direct factor $C_e$, namely the image of the subgroup $M_2(q)$ of $P\Gamma L_2(q)$ generated by $L_2(q)$ and an element $t\mapsto at^{\sqrt q}$ where $a$ is a non-residue. (For instance, when $q=9$ this is the `Mathieu group' $M_{10}=A_6.2_3$ in ATLAS notation~\cite{ATLAS}.) In this case there are $3$-transitive groups $G$ which do not contain $PGL_2(q)$; we need to show that none of them contains a cycle with two fixed points. We can take these to be $0$ and $\infty$, in which case the subgroup of $P\Gamma L_2(q)$ fixing them consists of the semilinear transformations $g:t\mapsto at^{\gamma}$, where $a\in{\mathbb F}_q^*$ and $\gamma\in{\rm Gal}\,{\mathbb F}_q$. This has a normal subgroup $N=\{g\mid\gamma=1\}\cong {\mathbb F}_q^*\cong C_{q-1}$, complemented by a subgroup $\{g\mid a=1\}\cong{\rm Gal}\,{\mathbb F}_q\cong C_e$. Replacing $g$ with a suitable power of the same order, we may assume that $\gamma:t\mapsto t^{p^f}$ for some $f$ dividing $e$, so $\gamma$ has order $d=e/f$. Then
\[g^d:t\mapsto a^{1+p^f+p^{2f}+\cdots+p^{(d-1)f}}t\]
is an element of $N$, and since $1+p^f+p^{2f}+\cdots+p^{(d-1)f}$ divides $p^e-1=q-1$ the order $m$ of $g^d$ divides
\[\frac{q-1}{1+p^f+p^{2f}+\cdots+p^{(d-1)f}}.\]
If $g$ is a cycle of length $q-1$ then its order $dm$ is equal to $q-1$, so
$1+p^f+p^{2f}+\cdots+p^{(d-1)f}$ divides $d$. This is impossible unless $d=1$, so that $g\in PGL_2(q)$. It follows that the only $(q-1)$-cycles $g\in P\Gamma L_2(q)$ are those in $PGL_2(q)$; since they satisfy $\langle L_2(q), g\rangle=PGL_2(q)$, the only groups $G$ containing $(q-1)$-cycles are those containing $PGL_2(q)$.

\section{Jordan and Marggraff}

Following Burnside~\cite[\S 159]{Bur} and Wielandt~\cite[Theorem 13.8]{Wie}, the result that a primitive permutation group containing a cycle with $k$ fixed points must be $(k+1)$-transitive has often been attributed to Marggraff (or Marggraf or Marggraaf). Both authors state the result without proof, referring to his dissertation~\cite{Mar}. Dixon and Mortimer~\cite[Exercise~7.4.11]{DM} set it as an exercise, without attribution or solution, though in a later hint they refer to a proof by Levingston and Taylor~\cite{LT}. In his scholarly review of that paper, Neumann~\cite{Neu} points out that an earlier paper by Atkinson~\cite{Atk} contains a similar proof due to Alan Williamson.

In fact Neumann, clearly one of the few who have read Marggraff's dissertation or his subsequent paper~\cite{Mar95}, argues {in~\cite{Neu} and in more detail in~\cite{Neu85} that this theorem should really be attributed to Jordan. Here is Jordan's Th\'eor\`eme~I from p.~384 of~\cite{Jor1} (see also~\cite[p.~314]{Jor2}), with the incorrect `{\it $n-p-2q+3$ fois primitif}' in his first sentence amended to `{\it $n-p-2q+3$ fois transitif}', the phrase he surely intended (see~\cite[p.~272]{Neu85} for Neumann's comments on this, including an English translation of  Th\'eor\`eme~I using modern terminology):

{\it Si un groupe $G$, primitif et de degr\'e $n$, contient un groupe $\Gamma$ dont les substitutions ne d\'eplacent que $p$ lettres et les permutent transitivement ($p$ \'etant un entier quelconque), il sera au moins $n-p-2q+3$ fois transitif, $q$ \'etant le plus grand diviseur de $p$ tel, que l'on puisse r\'epartir les lettres de $\Gamma$ de deux mani\`eres diff\'erentes en syst\`emes de $q$ lettres jouissant de la propri\'et\'e que chaque substitution de $\Gamma$ remplace les lettres de chaque syst\`eme par celles d'un meme syst\`eme.

Si aucun des diviseurs de $p$ ne jouit de cette propri\'et\'e (ce qui arrivera notamment si $\Gamma$ est primitif, ou form\'e des puissances d'une seule substitution circulaire), $G$ sera $n-p+1$ fois transitif.}

Note in particular the last sentence, which includes the case where $\Gamma$ is generated by a cycle. Neumann also finds no clear justification for the date of 1892 assigned by Burnside and Wielandt to Marggraff's dissertation, arguing that the rather sketchy evidence available suggests that it was probably written in 1889 or 1890. Again, see~\cite{Neu85} for more on Marggraff's work and its relationship with that of Jordan.

Concerning Jordan's Theorem~1.1, although Wielandt~\cite[Theorem 13.9]{Wie} refers to~\cite{Jor73}, it is not explicitly stated there. However, it follows easily from Th\'eor\`eme~I of~\cite{Jor1}, stated above, together with Th\'eor\`eme~I of~\cite{Jor73}:

{\it Soit $p$ un nombre premier impair. Un groupe de degr\'e $p+k$ ne pourra \^etre plus de $k$ fois transitif, si $k>2$, \`a moins de contenir le groupe altern\'e.}

\bigskip

\noindent{\bf Acknowledgements} The author thanks Alexander Zvonkin for raising the issue discussed here, the organisers of the conference Groups and Riemann Surfaces, Madrid, September 2012, where this interaction took place, and Peter Neumann for a number of very helpful mathematical, stylistic and historical comments on an early draft of this paper.}

\end{document}